\documentclass[twocolumn,letterpaper]{article}

\usepackage[margin=0.75in]{geometry} 
\usepackage{times}      
\usepackage{helvet}     
\usepackage{courier}    

\usepackage[hyphens]{url}
\usepackage{graphicx}
\usepackage{caption}
\usepackage{amsmath,amssymb}
\usepackage{booktabs}
\usepackage{multirow}
\usepackage{subcaption}
\usepackage{algorithm}
\usepackage{algpseudocode}

\usepackage{natbib}

\usepackage{authblk} 

\title{\bfseries Vibroseis Vehicle Routing Problem with Spatio-Temporal Coupled Constraints}

\author[1]{Kexin Zhu}
\author[1]{Jialong Shi\thanks{Corresponding author: jialong.shi@xjtu.edu.cn}}
\author[1]{Jianyong Sun}
\author[2]{Heng Zhou}
\author[1]{Mingen Kuang}
\author[3]{Ye Fan}

\affil[1]{School of Mathematics and Statistics, Xi'an Jiaotong University}
\affil[2]{BGP Inc., China National Petroleum Corporation}
\affil[3]{School of Electronics and Information, Northwestern Polytechnical University}

\date{}

\usepackage{bibentry}

\begin{document}
\maketitle

\begin{abstract}
Simultaneous multi vibroseis vehicle operations are central to modern land seismic exploration and can be modeled as a Vehicle Routing Problem (VRP). A critical distinction from classical VRPs, however, is the need for a minimum start-time interval between nearby sources to prevent signal interference. This constraint introduces strong spatio-temporal coupling, as one vehicle's route directly impacts the schedules of others, leading to potential forced waits. This paper defines this novel problem as the Vibroseis Vehicle Routing Problem with Spatio-Temporal Coupled Constraints (STCVRP) and aims to minimize the makespan.

To systematically investigate this problem, we first establish a Mixed-Integer Linear Programming (MILP) model to provide a precise mathematical description. As the MILP is intractable for realistic-scale instances, we subsequently develop a discrete-event simulation model. This model accurately calculates the true makespan of any candidate solution under all constraints, thereby serving as a high-fidelity fitness evaluation function for metaheuristics. To solve the problem, we design and implement a simulation-based Genetic Algorithm (GA) framework capable of searching the solution space efficiently. Finally, to validate the proposed framework and facilitate future research, we construct and release a new benchmark suite for the STCVRP. We solve the instances using both exact methods and our GA, with experimental results demonstrating the problem's high complexity and the effectiveness of the metaheuristic approach, providing a set of high-quality initial solutions for the benchmark.
\end{abstract}

\section{Introduction}
Seismic exploration is a method of acquiring subsurface stratigraphic morphology and physical properties by artificially generating seismic waves and analyzing their propagation through the earth. Vibroseis is a mechanical seismic source that generates seismic waves by continuously impacting the ground with vibrators mounted on vibroseis vehicles. Compared to traditional methods using explosives and heavy hammers, it offers advantages such as flexibility, safety, and environmental friendliness. Since Conoco Oil Company first applied vibroseis technology to seismic exploration in the 1950s \citep{dean2016high,laing1972some}, this technology has continuously improved in signal quality and operational efficiency, gradually becoming the mainstream seismic method.

In a typical vibroseis survey, a set of locations, known as task points, are first established across the exploration area. Vibroseis vehicles then traverse these points, stopping at each one to generate controlled seismic signals (namely, a sweep). Simultaneously, an array of sensors deployed in the field records the resulting seismic waves. An acquisition survey is considered complete once all designated task points have been serviced. To enhance acquisition efficiency, the slip-sweep technique\citep{rozemond1996slip} was proposed, which has been continuously refined in practical applications. Slip-sweep allows for continuous, overlapping sweeps, where multiple vehicle groups can work simultaneously with their waveforms separated by computer processing, significantly reducing time intervals between adjacent sweeps and notably improving the efficiency of vibroseis survey.\citep{chengfu2014implementation,jiaotong2018time}.

However, in such scenarios where multiple vibroseis vehicles work simultaneously, mutual harmonic interference can occur when vehicles are geographically proximate and perform tasks with temporal overlap, thereby affecting task execution or leading to failure \citep{bouska2010distance,bagaini2006overview}. Therefore, a minimum "separation" time must be enforced between the execution windows of any two tasks; such a constraint is called the \textit{slip-time separation rules}. This principle originated from the classical Vibroseis Slip-Sweep technique, whose core idea is to set a minimum time delay, i.e. "slip-time" \citep{xinquan2014new} for adjacent excited source groups to ensure their signals can be effectively distinguished and separated during recording. This technique has been widely applied in various advanced acquisition techniques, such as DSSS with slip-sweep and managed spread and source (MSS), where the required minimum slip time is typically a function of the physical distance between seismic sources \citep{dean2024practical,zhang2025high}.

The operational paradigm in vibroseis exploration, where multiple vehicles are dispatched to service their respective task points, can be modeled as a Vehicle Routing Problem (VRP). However, it is the unique spatio-temporal coupled constraint that makes this problem distinct and highly complex.The traditional VRP focuses on path optimization and capacity constraints, with numerous variants including the VRP with time windows and the multi-depot VRP \citep{toth2014vehicle,solomon1987algorithms,vidal2013heuristics,kitjacharoenchai2019multiple}. In recent years, VRP with spatio-temporal constraints have gained increasing attention, such as time-dependent VRP and multi-vehicle synchronization constrained VRP \citep{ichoua2003vehicle,gendreau2015time,drexl2012synchronization}. However, existing spatiotemporal constrained VRP research primarily concentrates on synchronization constraints that promote vehicle collaboration, such as requiring multiple vehicles to arrive at the same location simultaneously to perform tasks, without addressing time separation constraints based on physical distance.

The Spatio-Temporal Coupled Vibroseis Routing Problem (STCVRP) studied in this paper fundamentally differs from the existing VRP variants. Traditional synchronized VRPs aim to achieve coordinated cooperation among vehicles, whereas STCVRP requires ensuring spatio-temporal separation between vehicles to avoid harmonic interference. Specifically, STCVRP necessitates the simultaneous optimization of vehicles' physical paths and the precise execution schedule for all tasks, where the start time interval for vehicles at each task point must satisfy the function constraint defined by the slip-time separation rules. This physical distance-based time interval constraint makes STCVRP a highly complex optimization problem with strong spatio-temporal coupling.

The main contributions of this study are as follows: Firstly, we propose the Spatio-Temporal Coupled Vibroseis Vehicle Routing Problem (STCVRP), which introduces a new type of constraint to the VRP and fills a gap in vibroseis vehicle scheduling optimization research. Secondly, we establish a Mixed-Integer Linear Programming (MILP) model and a Simulation-Based Solution Evaluation Model for STCVRP, and develop a genetic algorithm-based solution method to address the high complexity of STCVRP. Furthermore, we construct a STCVRP benchmark dataset and make it publicly available, which aims to provide a standardized evaluation platform for further research on this problem.

The remainder of this paper is organized as follows: Section \uppercase\expandafter{\romannumeral 2} defines a general STCVRP problem, establishes an MILP mathematical model for its precise description, and proposes a Simulation-Based Solution Evaluation Model for subsequent solving. Section \uppercase\expandafter{\romannumeral 3} introduces the heuristic solution method for STCVRP. Section \uppercase\expandafter{\romannumeral 4} describes the method for generating STCVRP instances from existing benchmark test instances. Section \uppercase\expandafter{\romannumeral 5} presents the experimental results and analysis. Section \uppercase\expandafter{\romannumeral 6} provides the conclusions of this paper and outlines future research directions.

\section{The Spatio-Temporal Coupled Vibroseis Vehicle Routing Problem (STCVRP)}
\subsection{Problem Description}
The STCVRP is an NP-hard problem that can be described as follows: Given a set of vibroseis vehicles and a series of geographically dispersed task points, where each task point requires a vehicle to remain for a fixed duration to perform work, the objective is to plan the optimal path and precise service timetable for each vehicle while minimizing the maximum completion time and strictly adhering to spatio-temporal coordination constraints.

Unlike conventional VRPs, which primarily focus on minimizing cost or total travel time, the core challenge of STCVRP lies in simultaneously optimizing both the vehicles' physical paths and the precise execution schedule for all tasks. To optimize the excitation efficiency of vibroseis vehicles, it is crucial to minimize their total operational time. The objective function for STCVRP is to minimize the maximum vehicle completion time, which consists of three main components:
{\fontsize{10}{11}\selectfont
\begin{equation}
T = \max_{i} \left( T_{i}^{s} + T_{i}^{w} + T_{i}^{m} \right)
\end{equation}
}
where $T_i^s$ is the sweeping(working) time of the i-th vibroseis vehicle, $T_i^w$ is the waiting time incurred by the i-th vibroseis vehicle to avoid harmonic interference, and $T_i^m$ is the relocation time for the i-th vibroseis vehicle to move from its current task point to the next. 

As the working duration $T_i^w$ is typically a fixed value, the primary focus of optimization lies in minimizing the waiting time $T_i^r$ and the relocation time $T_i^m$.

\subsubsection{Problem Constraints}
The STCVRP is subject to two key constraints that distinguish it from traditional VRP variants:

\paragraph{(1) Fixed Service Time Constraint} 
Due to the nature of seismic wave excitation and scanning, each vibroseis vehicle must remain at a task point for a fixed duration to complete its work. We assume this service time is a known, fixed constant for all task points.

\begin{figure*}[t]
\centering
\includegraphics[width=2.0\columnwidth]{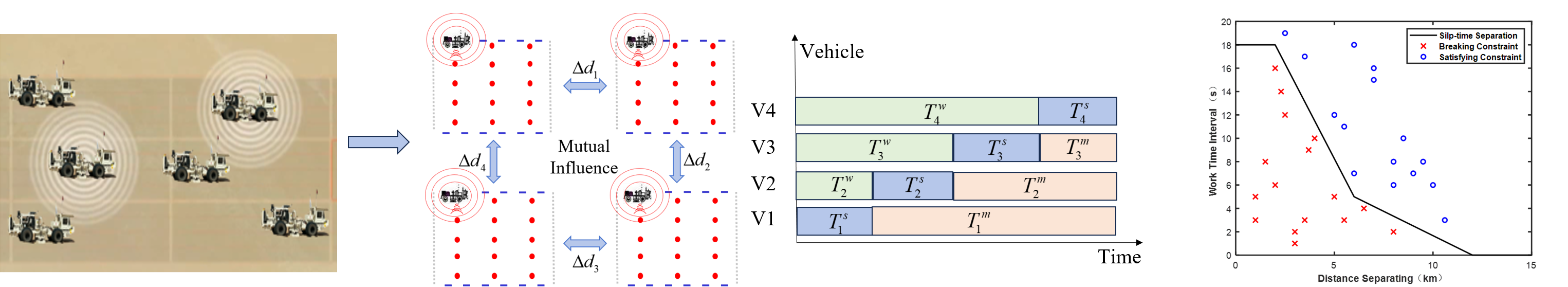} 
 \caption{Illustration of the slip-time separation rule. To avoid interference from Vehicle $v_1$, Vehicles $v_2$, $v_3$, and $v_4$ must wait and start their operations sequentially according to their priority order. And the left plot shows the rule function (black line) defining the minimum required start-time interval based on distance, with feasible (blue) and infeasible (red) vehicle pairs.}
\label{tab:silp-time rule}
\end{figure*}

\paragraph{(2) Slip-time Separation Constraint} 
To prevent harmonic interference, the start time difference between any two tasks performed by different vibroseis vehicles must satisfy the \textit{slip-time separation rules}. The required start time difference is solely dependent on the physical distance separating the task points. The separation rules can be represented as a piecewise function, as shown in Figure\ref{tab:silp-time rule}. Only points above the rule function line represent valid operating conditions that satisfy the constraint.

A crucial mechanism to ensure adherence to these rules is the introduction of waiting time. If the planned start times of two vehicles violate the separation rules, one vehicle can commence work while the other(s) wait until the minimum required time interval is satisfied. By incorporating waiting time, virtually any task assignment and routing scheme can be transformed into a feasible schedule.The schematic diagram of this problem is shown in \ref{tab:silp-time rule}.
\subsubsection{Problem Assumptions}
This study operates under a set of standard VRP assumptions combined with problem-specific simplifications:

 all vehicle tours start and end at a central depot, and each task point is visited exactly once. We also assume unlimited vehicle resources (i.e., no fuel or mileage constraints), a constant travel speed, and that waiting is only permitted at task points, not en route. All vehicles depart from the depot at time zero.

\subsection{MILP Problem Formulation}
The STCVRP is defined on a complete directed graph \(G=(N_0, E)\). The node set \(N_0 = C \cup \{0\}\) comprises \(N\) task points \(C = \{1, 2, \dots, N\}\) and a central depot (node 0). The arc set \(E = \{(i, j) \mid i, j \in N_0, i \neq j\}\) contains all potential travel paths between nodes.

Let \(K_{max}\) be the number of available vehicles, which are assumed to be homogeneous and must all be utilized for service. The STCVRP mandates that when any two task points \(i\) and \(j\) are served by different vehicles, their respective service start times must maintain a minimum time interval \(g_{ij}\). The model employs binary variables \(x_{ijk}\) (1 if vehicle \(k\) traverses arc \((i, j)\), 0 otherwise) and \(v_{ik}\) (1 if task point \(i \in C\) is assigned to vehicle \(k\), 0 otherwise). An auxiliary binary variable \(z_{ij}\) indicates whether task points \(i\) and \(j\) are served by the same vehicle. Continuous variables include arrival time \(b_i\), waiting time \(t_i\), and service start time \(s_i\) at task point \(i \in C\), with \(s_i = b_i + t_i\). The completion time for vehicle \(k\) is denoted as \(C_k\), and the makespan as \(T\). The notation for the STCVRP model is summarized in Table \ref{tab:notation}.

\begin{table}[t]
\centering
\setlength{\tabcolsep}{1mm}
\small
\begin{tabular}{@{}ll@{}}
\toprule
\multicolumn{2}{l}{Indices and Sets} \\
\midrule
\(C\)           & Set of task points, \(C = \{1, 2, \dots, N\}\) \\
\(K\)           & Set of available vehicles, \(K = \{1, 2, \dots, K_{max}\}\) \\
\(N_0\)         & Set of all nodes including the depot, \(N_0 = C \cup \{0\}\) \\
\(i, j, h\)     & Indices for nodes \\
\(k\)           & Index for vehicles \\
\midrule
\multicolumn{2}{l}{Parameters} \\
\midrule
\(N\)           & Total number of task points \\
\(K_{max}\)     & Number of available vehicles \\
\(c_{ij}\)      & Travel time from node \(i\) to node \(j\) \\
\(w\)           & Uniform service time at each task point \\
\(g_{ij}\)      & Minimum time separation for services at \(i, j\) \\
& by different vehicles \\
\(M\)           & A sufficiently large positive constant \\
\midrule
\multicolumn{2}{l}{Decision Variables} \\
\midrule
\(x_{ijk}\)     & Binary, 1 if vehicle \(k\) travels from node \(i\) to node \(j\) \\
\(v_{ik}\)      & Binary, 1 if task point \(i\) is served by vehicle \(k\) \\
\(u_k\)         & Binary, 1 if vehicle \(k\) is used \\
\(z_{ij}\)      & Binary, 1 if task points \(i\) and \(j\) are served by\\
                & the same vehicle \\
\(u'_{ijk}\)    & Auxiliary non-negative variable for linearization \\
\(b_i\)         & Continuous, arrival time at task point \(i\) \\
\(t_i\)         & Continuous, waiting time at task point \(i\) \\
\(s_i\)         & Continuous, service start time at task point \(i\) \\
\(C_k\)         & Continuous, completion time of vehicle \(k\) \\
\(T\)           & Continuous, the makespan \\
\bottomrule
\end{tabular}
\caption{Notation for the STCVRP Model.}\label{tab:notation}
\end{table}

The MILP formulation for STCVRP is as follows:

\paragraph{Objective Function:}
{\fontsize{10}{11}\selectfont
\begin{equation}
\min T
\end{equation}
}

\paragraph{Routing and Assignment Constraints:}
{\fontsize{10}{11}\selectfont
\begin{align}
    & \sum_{k \in K} v_{ik} = 1 && \forall i \in C \\
    & v_{ik} = \sum_{j \in N_0, j \neq i} x_{ijk} && \forall i \in C, \forall k \in K \\
    & v_{ik} = \sum_{j \in N_0, j \neq i} x_{jik} && \forall i \in C, \forall k \in K \\
    & \sum_{i \in N_0, i \neq h} x_{ihk} - \sum_{j \in N_0, j \neq h} x_{hjk} = 0 && \forall h \in C, \forall k \in K \\
    & u_k = \sum_{j \in C} x_{0jk} && \forall k \in K \\
    & u_k = 1 && \forall k \in K
\end{align}
}

\paragraph{Time Flow and Subtour Elimination Constraints}
{\fontsize{10}{11}\selectfont
\begin{align}
    & s_i = b_i + t_i && \forall i \in C \\
    & b_j \ge s_i + w + c_{ij} - M \left(1 - \sum_{k \in K} x_{ijk}\right) && \forall i, j \in C, i \neq j \\
    & b_j \ge c_{0j} - M \left(1 - \sum_{k \in K} x_{0jk}\right) && \forall j \in C
\end{align}
}

\paragraph{Inter-Vehicle Interference Constraints:}
{\fontsize{10}{11}\selectfont
\begin{align}
    & u'_{ijk} \ge v_{ik} - v_{jk} && \forall i, j \in C (i<j), \forall k \in K \\
    & u'_{ijk} \ge v_{jk} - v_{ik} && \forall i, j \in C (i<j), \forall k \in K \\
    & 2 \cdot (1 - z_{ij}) = \sum_{k \in K} u'_{ijk} && \forall i, j \in C (i<j) \\
    & s_j - s_i \ge g_{ij} - M \cdot z_{ij} && \forall i, j \in C (i < j) \\
    & s_i - s_j \ge g_{ij} - M \cdot z_{ij} && \forall i, j \in C (i < j)
\end{align}
}

\paragraph{Makespan Constraints:}
{\fontsize{10}{11}\selectfont
\begin{align}
    & C_k \ge s_i + w + c_{i0} - M(1-x_{i0k}) && \forall i \in C, \forall k \in K \\
    & T \ge C_k && \forall k \in K
\end{align}
}
The objective function minimizes the makespan \(T\), encouraging both efficient routing and load balancing among vehicles. Constraints (4)-(9) ensure each task point is served by exactly one vehicle while maintaining routing consistency and flow conservation. All \(K_{max}\) vehicles must be utilized as specified by constraint (9). Constraints (10)-(12) manage time propagation along routes and eliminate subtours through the monotonic property of time variables. The introduction of waiting time \(t_i\) provides the necessary flexibility for the model to handle subsequent scheduling conflicts. Constraint (11) ensures the correct propagation of time along a path and effectively eliminates all sub-tours that do not include the depot through the monotonic increasing property of time variables.

The core innovation lies in constraints (13)-(17), which handle inter-vehicle interference. Constraints (13)-(15) collectively establish a rigorous logical relationship to determine the value of the auxiliary variable \(z_{ij}\). This is equivalent to \(2 \cdot (1-z_{ij}) = \sum_{k \in K} |v_{ik} - v_{jk}|\), and it has been linearized by introducing the auxiliary variable \(u'_{ijk}\). The effect is that \(z_{ij}\) is 1 if and only if task points \(i\) and \(j\) are served by the same vehicle; if they are served by different vehicles, \(z_{ij}\) is 0. Subsequently, constraints (16) and (17) utilize \(z_{ij}\) as a switch to activate the interference conditions. If \(z_{ij}=1\) (same-vehicle service), the \(M \cdot z_{ij}\) term renders these two constraints relaxed and inactive, which is logical as a vehicle will not experience harmonic interference with itself. However, if \(z_{ij}=0\) (different-vehicle service), these two constraints are activated, jointly enforcing \(|s_i - s_j| \ge g_{ij}\).

Constraints (18) and (19) are used to calculate and minimize the final objective value. Constraint (18) identifies the last task point on the path of each vehicle and calculates the total completion time \(C_k\) for that vehicle based on its service completion time and return travel time. Finally, constraint (19) defines the objective variable as the maximum among all vehicle completion times \(C_k\).

\subsection{Simulation-Based Solution Evaluation Model}
The computational complexity of the MILP model for the STCVRP is prohibitively high, rendering it intractable for the large-scale instances typical of real-world operational planning. We found that under strong interference conditions, obtaining the optimal solution for problems with just eight task points took over an hour, and larger instances could not be solved within a reasonable timeframe. A key reason for this difficulty is that the objective function of the STCVRP is \textbf{not in a closed-form}; the total waiting time for a given route set is a result of the system's dynamic interactions, not a simple sum of pre-calculable components. This characteristic necessitates a simulation-based approach for accurate evaluation.

To improve efficiency over simpler time-stepping methods, we adopted a Discrete-Event Simulation (DES) paradigm to develop the Event-Driven Scheduler with Spatiotemporal Constraints (EDSSC)\citep{dauzere2024flexible,burdett2018integrated}. This scheduler is a purpose-built evaluation model for the STCVRP, designed specifically to resolve the \textbf{cascading delay effects} inherent to the problem, where a single waiting decision can propagate through the entire multi-vehicle system. The model's architecture is event-driven, managing a priority queue of three key event types: \texttt{ARRIVE}, \texttt{START\_WORK}, and \texttt{END\_WORK}. The complete behavior of the EDSSC is formally defined by Algorithm \ref{alg:scheduler} and Algorithm \ref{alg:arrive_handler}, which details the top-level process, and Algorithm \ref{alg:conflict_resolution}, which specifies the core conflict resolution logic.
\begin{algorithm}[t]
\caption{Event-Driven Scheduler with Spatiotemporal Constraints}
\label{alg:scheduler}
\textbf{Input}: A solution $S$, problem instance $I$ \\
\textbf{Output}: The makespan for solution $S$
\begin{algorithmic}[1]
    \State $Q \gets \text{Initialize event queue with initial ARRIVE events}$\\$\text{ for all vehicles}$
    \While{$Q$ is not empty}
        \State $e \gets Q.pop()$
        \If{$e.type = \text{ARRIVE}$}
            \State \Call{HandleArriveEvent}{$e, Q, I$}
        \ElsIf{$e.type = \text{START\_WORK}$}
            \State $v \gets e.\text{vehicle}$; \quad $v.isWorking \gets \text{true}$
            \State $Q.push(\text{WorkEvent}(e.time + I.serviceTime))$
        \ElsIf{$e.type = \text{END\_WORK}$}
            \State $v \gets e.\text{vehicle}$; \quad $v.isWorking \gets \text{false}$
            \If{$v$ has a next task}
                \State $Q.push(\text{ArriveEvent}(e.time))$
            \Else
                \State Record completion time for vehicle $v$
            \EndIf
        \EndIf
    \EndWhile
    \State \textbf{return} max completion time over all vehicles
\end{algorithmic}
\end{algorithm}

\begin{algorithm}[t]
\caption{HandleArriveEvent}
\label{alg:arrive_handler}
\textbf{Input}: Current arrival event $e_{current}$, event queue $Q$, instance $I$ \\
\textbf{Output}: Updated $Q$ with scheduled START\_WORK events
\begin{algorithmic}[1]
    \State $batchArrives \gets \begin{tabular}[t]{@{}l@{}} 
                               \text{Extract concurrent ARRIVE events} \\ 
                               \text{from } Q 
                           \end{tabular}$
    \State Sort $batchArrives$ by priority (tasks completed asc, then vehicle ID asc)
    \For{each event $e_i$ in sorted $batchArrives$}
        \State $v_i \gets e_i.\text{vehicle}$
        \State $s_i \gets$ \Call{CalculateStartTime}{$v_i, e_i.time, I$} 
        \State Update $v_i$'s state: $isWorking \gets \text{true}$, $expectedStartTime \gets s_i$, etc.
        \State $Q.push(\text{CreateStartWorkEvent}(s_i))$
    \EndFor
\end{algorithmic}
\end{algorithm}

\begin{algorithm}[t]
\caption{CalculateStartTime(Conflict Resolution)}
\label{alg:conflict_resolution}
\textbf{Input}: Current vehicle $v_i$, arrival time $t_{arrival}$, instance $I$ \\
\textbf{Output}: The earliest valid start time $s_i^*$ for $v_i$
\begin{algorithmic}[1]
   \State $s_i^{\text{cand}} \gets t_{arrival}$
    \Repeat
        \State $s_i^{\text{prev}} \gets s_i^{\text{cand}}$
        \For{each other active vehicle $v_j$}
            \State $s_j \gets v_j.\text{expectedStartTime}$; \quad $e_j \gets v_j.\text{expectedEndTime}$
            \If{spatiotemporal conflict exists between $v_i$ (at $s_i^{\text{cand}}$) and $v_j$}
                \State $g_{ij} \gets \text{RequiredSeparationTime}(v_i, v_j)$
                \State $s_i^{\text{cand}} \gets \max(s_i^{\text{cand}}, \min(s_j + g_{ij}, e_j))$
            \EndIf
        \EndFor
    \Until{$s_i^{\text{cand}}$ converges}
    \State \textbf{return} $s_i^{\text{cand}}$
\end{algorithmic}
\end{algorithm}
The EDSSC employs two key mechanisms to handle the complexities of the STCVRP. First, a concurrent event handling mechanism (Algorithm \ref{alg:scheduler}, \ref{alg:arrive_handler}) addresses situations where multiple vehicles arrive simultaneously. Its prioritized batch processing is critical for ensuring deterministic and consistent evaluation, as an arbitrary handling of simultaneous events could lead to different makespans for the exact same solution, rendering the optimization search unreliable. Second, the core iterative conflict resolution algorithm (Algorithm \ref{alg:conflict_resolution}) calculates the earliest legitimate start time for any given vehicle. This iterative refinement is essential for the model's accuracy; a non-iterative, single-pass check would fail to capture secondary conflicts, systematically underestimating the true makespan and providing misleading guidance to the heuristic search.

This simulation model will serve as the fitness evaluation function for the Genetic Algorithm described subsequently, providing accurate guidance for its evolutionary search.

\section{Optimization Method}
To optumize the STCVRP, we adopt a simulation-based genetic algorithm framework, utilizing the simulation model (EDSSC) described in Section 1.3 as the core fitness evaluation method. Due to the strong spatio-temporal coupling, changes in a single vehicle's tasks can create cascading effects that impact the overall completion time. The genetic algorithm employs the simulation scheduler as a black-box fitness function, enabling accurate evaluation of complex spatio-temporal interactions that cannot be captured by traditional distance-based fitness calculations, to search for vehicle routing plans that minimize the overall completion time within the complex solution space.

\subsubsection{Encoding and Fitness Evaluation}
Each chromosome is represented as a two-dimensional list, where each sublist represents a vehicle's visit sequence of task points. For example, [[t1, t2], [t3, t4, t5]] indicates that Vehicle 1's path is depot → t1 → t2 → depot, and Vehicle 2's path is depot → t3 → t4 → t5 → depot. The fitness value of an individual is defined as the minimum completion time calculated by the EDSSC simulation, with smaller fitness values indicating higher solution quality.

\subsubsection{Genetic Algorithm Design}
\textbf{\textit{Initialization Strategy:}} To ensure the diversity of initial solutions, we adopt multiple initialization strategies: K-means based spatial clustering initialization, nearest neighbor greedy construction, balanced load allocation, and random generation.
\textbf{\textit{Genetic Operations:}} Tournament selection is used for parent selection. The crossover operation uses the classic Order Crossover (OX1) operator \citep{davis1985applying,arnold2019efficiently} to preserve good route structures from parents. The mutation operation adopts a hybrid strategy, where each mutation randomly selects either 2-opt local search \citep{croes1958method,nakamura2020short} for route optimization or optimal position insertion mutation for global search. Notably, the insertion mutation uses approximate fitness evaluation based on Euclidean distance to avoid costly simulation calls, significantly improving mutation efficiency. Population update employs an elitism strategy.

The overall algorithm flow is shown in Appendix. By using the EDSSC simulator as an accurate fitness evaluation function, this genetic algorithm framework can efficiently search for high-quality STCVRP solutions that satisfy spatio-temporal separation constraints.

\section{STCVRP Benchmark Test Suite}
To validate the problem's complexity and provide a benchmark for future research, we designed a dedicated STCVRP dataset. It employs two strategies: modifying classic VRP datasets (Solomon and 'rat' from TSPLIB \citep{solomon1987algorithms,reinelt1991tsplib,gunawan2021vehicle}) to leverage their node distribution characteristics, and generating grid instances with random perturbations to simulate common grid-like layouts in exploration.

\subsection{Parameter Configuration}
We establish fixed parameters to ensure standardized testing conditions. The target average task point distance is 40 meters, based on nearest neighbor distances in operational areas. Vehicle speed is 5.0 m/s. Both service time and maximum waiting time are set to 8 seconds, forming the core spatio-temporal constraints. The constraint intensity is controlled by adjusting the maximum constraint distance.

\subsection{Instance Generation}
For transformed instances, we extract coordinates from classic datasets, calculate the average nearest neighbor distance, and apply a scaling factor of $40.0 / \text{current\_avg\_distance}$ to all coordinates while preserving the distribution. For grid-based instances, we generate a near-square grid based on the number of customers and target distance, then add Gaussian noise to each node to mimic real-world deviations.

\subsection{Benchmark Instances}
The constructed test dataset contains instances at two scales under uniform parameter setup. Small-scale instances $(N<100)$ contain three distribution patterns: clustered (C), random (R), and grid (G), with coordinates normalized to an average nearest-neighbor distance of 40 meters. Medium-scale instances utilize TSPLIB datasets 'rat575' and 'rat783', featuring grid distributions with node perturbations.

Parameter selection addresses algorithm performance and sensitivity evaluation. The vehicle-to-customer ratio ranges from $1:5$ to $1:12.5$ for small-scale instances and extends to $1:20$ to $1:40$ for medium-scale instances. The maximum constraint distance ($d$) determines constraint intensity: 80 meters represents weak constraints, 150-200 meters indicates medium-to-strong constraints, and 500-800 meters signifies strong constraints requiring coordination among nearly all vehicles\citep{uchoa2017new}.

Each instance follows the naming convention [Task Point Distribution][Number of tasks]$\_$[Number of Vehicles]k$\_$[Maximum Constraint Distance]d, where 'C', 'R', and 'G' denote clustered, random, and grid distributions respectively\citep{mavrovouniotis2020benchmark}. The minimum operational time interval uses a piecewise linear function:
{\fontsize{9}{10}\selectfont
\begin{equation}
W(d) = \begin{cases} 
    W_{max} - \dfrac{W_{max}}{D_{max}} \cdot d & \text{if } 0 \le d < D_{max} \\
    0 & \text{if } d \ge D_{max}
\end{cases}
\end{equation}
}
where $W(d)$ is the minimum time interval varying with distance $d$, and $D_{max}$ is the maximum constraint distance. The detail of our benchmark instances is shown in Appendix.

\section{Experimental Results}
\subsection{Experimental Setup}
The MILP experiments were performed under a Windows System with an Intel Core i7-12700K (3.60GHz) processor and 16GB of RAM, using the Gurobi 12.0.1 solver\citep{gurobi2024}. A computational time limit of 2000 seconds was set for each instance. The simulation-based Genetic Algorithm (GA) was implemented in C++ and executed on the Tianhe-2 supercomputer, with 10 independent runs performed for each instance. 

The GA parameters were: population size $N_{pop}$ (50 for small-scale, 100 for medium-scale), crossover rate $p_c=0.8$, mutation rate $p_m=0.2$, elite individuals $N_e$ (2 for small-scale, 5 for medium-scale), with termination after 1000 consecutive iterations without improvement.

\subsection{Analysis of the Results}
Since the STCVRP is a newly-defined problem, we evaluate the performance of our heuristic by comparing its solutions with those obtained by the MILP. The results for small-scale and medium-scale instances are presented in Table \ref{tab:small_scale_results} and \ref{tab:medium_scale_results}, respectively. Note that the MILP approach failed to provide a feasible solution for some instances (N=100) and could not be applied to medium-scale problems. For the GA, all results are based on 10 independent runs, and the GAP (\%) relative to the MILP makespan is calculated where applicable.
\begin{table*}[t]
    \centering
    \begin{tabular}{@{}lrr|rrrrrr@{}}
        \toprule
        \multirow{2}{*}{\textbf{Instance}} & \multicolumn{2}{c|}{\textbf{MILP}} & \multicolumn{6}{c}{\textbf{Simulation-based GA}} \\
        \cmidrule(lr){2-3} \cmidrule(lr){4-9}
        & \textbf{Makespan} & \textbf{$\text{T}_\text{W(s)}$} & \textbf{Best} & \textbf{Worst} & \textbf{Mean±std} & \textbf{$\text{t}_\text{avg(s)}$} & \textbf{$\text{T}_\text{W(s)}$} & \textbf{GAP (\%)} \\
        \midrule
        C25\_2k\_150d  & \textbf{362.81} & 7.0  & 387.6 & 442.2 & 417.07 $\pm$ 16.48 & 0.39 & 0.13  & 6.83 \\
        C25\_5k\_150d  & \textbf{264.16} & 33.0 & 274.3 & 309.2 & 281.03 $\pm$ 12.05 & 0.37 & 3.72  & 3.84 \\
        
        R25\_2k\_150d  & 275.52 & 7.0  & \textbf{233.9} & 273.4 & 266.78 $\pm$ 11.74 & 0.30 & 4.74  & -15.11 \\
        R25\_5k\_150d  & 120.89 & 23.0 & \textbf{116.0} & 131.0 & 127.73 $\pm$ 4.48  & 0.41 & 12.96 & -4.05 \\
        
        G25\_2k\_150d  & 259.35 & 2.0  &\textbf{ 238.7 }& 258.8 & 248.61 $\pm$ 5.61  & 0.32 & 3.18  & -7.96 \\
        G25\_5k\_150d  &\textbf{ 130.96 }& 31.0 & 140.7 & 145.7 & 142.62 $\pm$ 1.59  & 0.36 & 12.20 & 7.44 \\
        \midrule
        C50\_5k\_150d  & 498.05 & 38.0 & \textbf{360.0} & 514.2 & 392.46 $\pm$ 59.80 & 0.76 & 1.38  & -27.72 \\
        R50\_5k\_150d  & 293.45 & 37.0 & \textbf{263.6} & 308.1 & 292.07 $\pm$ 14.43 & 0.87 & 6.20  & -10.17 \\
        G50\_5k\_150d  & 275.21 & 33.0 & \textbf{245.7} & 276.5 & 267.55 $\pm$ 9.05  & 0.73 & 19.81 & -10.72 \\
        
        C100\_8k\_150d & -      & -    & 865.6 & 996.6 & 964.77 $\pm$ 48.14 & 1.86 & 19.49 & - \\
        R100\_8k\_150d & -      & -    & 397.3 & 415.3 & 408.25 $\pm$ 8.64  & 1.86 & 33.25 & - \\
        G100\_8k\_150d & -      & -    & 342.4 & 349.2 & 343.88 $\pm$ 2.01  & 1.75 & 47.91 & - \\
        \midrule
        G50\_5k\_80d   & 296.44 & 31.0 & \textbf{247.9} & 260.3 & 255.74 $\pm$ 2.96  & 0.76 & 5.24  & -16.38 \\
        G50\_5k\_200d  & 283.05 & 57.0 & \textbf{262.3} & 285.2 & 275.57 $\pm$ 6.87  & 0.85 & 40.06 & -7.33 \\
        G50\_8k\_150d  & 195.53 & 60.0 & \textbf{193.8} & 241.7 & 219.99 $\pm$ 15.73 & 0.87 & 31.27 & -0.88 \\
        \bottomrule
    \end{tabular}
    \caption{Comparison of results for small-scale instances. GAP (\%) = $\tfrac{\text{GA}_{\text{Best}} - \text{MILP}_{\text{Makespan}}}{\text{MILP}_{\text{Makespan}}} \times 100\%$. - indicates that the exact method failed to find a feasible solution within the time limit.}
    \label{tab:small_scale_results}
\end{table*}

\begin{table}[t]
    \centering
    \setlength{\tabcolsep}{1mm}
    \small
    \begin{tabular}{@{}lrrrrr@{}}
        \toprule
        \multirow{2}{*}{\textbf{Instance}} & \multicolumn{5}{c}{\textbf{Simulation-based GA}} \\
        \cmidrule(lr){2-6}
        & \textbf{Best} & \textbf{Worst} & \textbf{Mean±std} & \textbf{$\text{t}_\text{avg(s)}$} & \textbf{$\text{T}_\text{W(s)}$} \\
        \midrule
        G575\_8k\_500d  & 1734.4 & 1786.5 & 1761.63 $\pm$ 16.81 & 102.35 & 422.37 \\
        G575\_15k\_120d & 1081.7 & 1140.8 & 1101.79 $\pm$ 20.71 & 77.17  & 90.14 \\
        G575\_15k\_500d & 1334.9 & 1589.7 & 1424.52 $\pm$ 90.76 & 96.01  & 1224.84 \\
        G575\_15k\_800d & 1765.6 & 1915.3 & 1843.28 $\pm$ 49.40 & 203.01 & 3905.24 \\
        G575\_30k\_500d & 1058.8 & 1144.6 & 1098.23 $\pm$ 26.86 & 168.05 & 4037.07 \\
        R575\_15k\_500d & 1673.1 & 1774.6 & 1711.14 $\pm$ 31.75 & 103.66 & 1091.82 \\
        \midrule
        G783\_20k\_120d & 1130.6 & 1160.1 & 1144.97 $\pm$ 7.87  & 173.33 & 169.09 \\
        G783\_20k\_500d & 1528.7 & 1781.4 & 1618.40 $\pm$ 89.94 & 207.52 & 1726.03 \\
        G783\_20k\_800d & 2019.6 & 2138.0 & 2088.84 $\pm$ 38.28 & 330.36 & 6077.78 \\
        G783\_40k\_500d & 1169.6 & 1269.4 & 1228.91 $\pm$ 29.00 & 336.55 & 6110.85 \\
        R783\_20k\_500d & 1951.1 & 2085.2 & 2026.54 $\pm$ 39.85 & 220.85 & 1542.01 \\
        \bottomrule
    \end{tabular}
    \caption{Results for medium-scale instances. MILP omitted due to computational intractability within time limits.} 
    \label{tab:medium_scale_results}
\end{table}

\subsubsection{Exact Algorithm Complexity Analysis}
To gain deeper insights into the computational complexity of STCVRP, we conducted exact solution experiments on small-scale problems using the Gurobi 12.0.1 solver. We tested instances with 6, 8, and 10 task points, with inter-task distances set to 1 unit.

\begin{figure}[t]
\centering
\includegraphics[width=0.9\columnwidth]{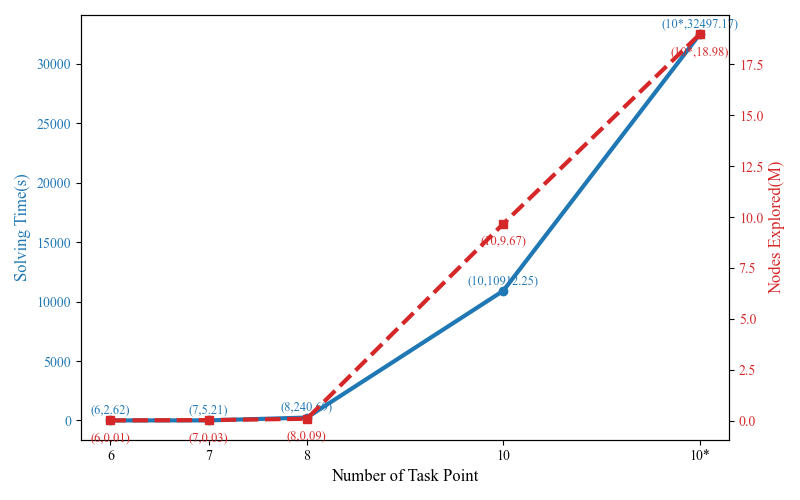} 
\caption{Computational performance of the MILP solver on extremely small-scale STCVRP instances.}
\label{tab:result}
\end{figure}

The experimental results demonstrate the exponential complexity of the problem: solving instances with 6 task points required 2.62 seconds, 8 task points required 240.69 seconds, while 10 task points required 10,912 seconds, showing exponential growth in solution time. The results shows in Figure \ref{tab:result}.More significantly, the constraint sensitivity is extreme: when the average distance between 10 task points was doubled (i.e., inter-task distance changed from 1 to 2), the solution time surged to 32,497 seconds, nearly tripling the computational effort, with the number of explored nodes increasing from 9.67 million to 18.98 million. These results indicate that STCVRP exhibits not only exponential scale complexity but also extreme sensitivity to constraint tightness, further validating the computational intractability of the problem and the necessity of heuristic algorithms.

\subsubsection{Heuristic Algorithm Performance Analysis}
An analysis of the experimental data reveals several key findings regarding the characteristics of the STCVRP and the performance of the proposed algorithm.

\begin{figure}[t]
\centering
\includegraphics[width=0.9\columnwidth]{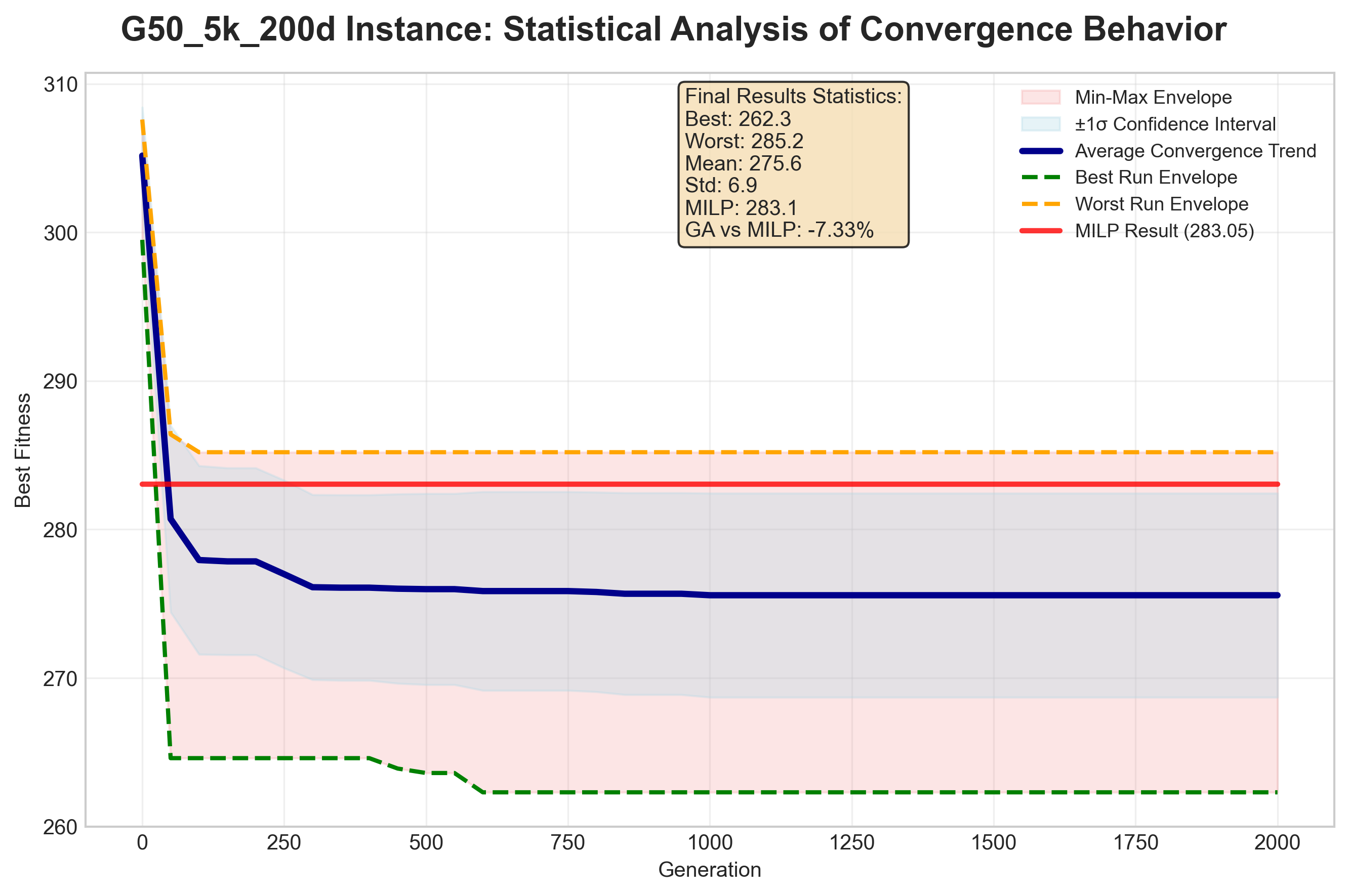} 
\caption{G50\_5k\_200d Instance: Individual Run Convergence Curves}
\label{tab:RunConvergenceCurves}
\end{figure}

\paragraph{Algorithm Performance and Computational Complexity}
The MILP approach could not find proven optimal solutions for instances with 25 tasks within the given time limit, and failed completely when the problem size increased to 100 tasks. The significant standard deviation observed in the GA's multiple runs indicates a complex solution landscape with numerous local optima. In contrast, our designed GA consistently finds high-quality solutions across all instances, requiring only 0.3-1.9 seconds for small-scale problems and demonstrating good scalability, highlighting the necessity of efficient metaheuristics for realistic-scale STCVRP instances. The convergence curves for the G50\_5k\_200d instance are shown in Figure \ref{tab:RunConvergenceCurves}.

\textbf{\textit{Impact of Spatio-temporal Constraints:}}
The spatio-temporal coupling constraint is a defining feature of the STCVRP, and its strength significantly affects the solution. This is evident in the medium-scale instances: as the maximum constraint distance increases from 120m to 800m, the total waiting time increases exponentially, with over 40-fold increases observed in multiple instance series. This quasi-exponential growth in waiting time directly leads to substantial increases in makespan, demonstrating that stronger spatio-temporal constraints result in increased inter-vehicle interference.

\textbf{\textit{System Design Trade-offs and Scaling Behavior:}}
The configuration of vehicle fleet size demonstrates clear trade-offs in the STCVRP. On medium-scale instances, increasing vehicles from 8 to 30 improves makespan due to improved parallelization, but also increases coordination overhead substantially. While doubling the fleet from 15 to 30 vehicles improved makespan by only about 20.7\%, the average waiting time per vehicle increased significantly, indicating diminishing marginal utility due to increased coordination overhead. The geographical distribution of task points also influences solution complexity: clustered distributions produce longer makespans and larger standard deviations due to intensive coordination requirements, while random distributions primarily challenge route optimization due to irregularity. The algorithm's computation time demonstrates super-linear scaling with problem size, increasing over 100-fold when the number of tasks increases from 50 to 575.

\section{Conclusions}
In this study, we introduced the Spatio-Temporal Coupled Vibroseis Routing Problem (STCVRP), a novel VRP variant characterized by physical distance-based time separation constraints. We established comprehensive mathematical models, including an MILP formulation and a discrete-event simulation model (EDSSC), and developed a simulation-driven genetic algorithm. Our experimental results demonstrate that STCVRP exhibits extremely high computational complexity and is highly sensitive to constraint tightness. The proposed heuristic algorithm proved capable of efficiently obtaining high-quality solutions across instances of various scales, validating its effectiveness for this challenging problem.

Future research can be extended in several directions. These include developing more efficient exact algorithms and hybrid heuristic methods, extending the spatio-temporal separation constraint model to other application domains.

\bibliographystyle{plainnat}
\bibliography{main}
\onecolumn
\begin{appendix}
\section{Introduction}
This appendix provides additional details and supplementary materials for the paper "Vibroseis Vehicle Routing Problem with Spatio-Temporal Coupled Constraints". The appendix is organized into four main sections.

The structure of this appendix is as follows: Appendix A presents the complete pseudocode of the genetic algorithm implementation, providing detailed algorithmic procedures that complement the high-level description in the main paper. Appendix B contains comprehensive information about the benchmark test instances, including detailed specifications and parameter configurations for both small-scale and medium-scale problem instances. Appendix C reports the detailed computational results of exact algorithm experiments on extremely small-scale instances. Finally, Appendix D illustrates the convergence behavior and performance characteristics of the genetic algorithm through iteration curves and statistical analysis.

These supplementary materials are intended to enhance the reproducibility of our research and provide researchers with comprehensive information for further investigation and comparison studies in the field of spatio-temporal constrained vehicle routing problems.

\section{Simulation-based Genetic Algorithm Framework}
The overall algorithm flow is shown in Algorithm \ref{alg:ssga}. By using the EDSSC simulator (explained in the main text) as an accurate fitness evaluation function, this genetic algorithm framework can efficiently search for high-quality STCVRP solutions that satisfy spatio-temporal separation constraints.

\begin{algorithm}[H]
\caption{Simulation-based GA Algorithm Framework}
\label{alg:ssga}
\begin{algorithmic}[1]
\Require Problem instance $I$, population size $N_{pop}$, max generations $G_{max}$
\Ensure The best found solution $S^*$ with the minimum makespan

\Procedure{GA}{$I, N_{pop}, G_{max}$}
    \State $P \gets \Call{InitializePopulation}{I, N_{pop}}$ \Comment{Create initial population}
    \For{each solution $S_i$ in $P$}
        \State $f_i \gets \Call{RunEDSSC}{S_i, I}$ \Comment{Using EDSSC to evaluate fitness}
    \EndFor
    
    \For{$g \gets 1$ to $G_{max}$}
        \State $P_{new} \gets \emptyset$
        \State $P_{elites} \gets \text{select best individuals from } P$ \Comment{Elitism}
        \State $P_{new} \gets P_{new} \cup P_{elites}$
        
        \While{$|P_{new}| < N_{pop}$}
            \State $(p_1, p_2) \gets \Call{SelectParents}{P}$ \Comment{e.g., Tournament Selection}
            \If{$\text{random}(0,1) < \text{crossoverRate}$}
                \State $(c_1, c_2) \gets \Call{Crossover}{p_1, p_2}$ \Comment{Order Crossover}
            \Else
                \State $(c_1, c_2) \gets (p_1, p_2)$
            \EndIf
            
            \State $\Call{Mutate}{c_1}$; \quad $\Call{Mutate}{c_2}$ \Comment{Hybrid Mutation}
            \State Add $c_1, c_2$ to $P_{new}$
        \EndWhile
        
        \State $P \gets P_{new}$
        \For{each solution $S_i$ in $P$}
            \State $f_i \gets \Call{RunEDSSC}{S_i, I}$ \Comment{Evaluate new population}
        \EndFor

        \If{termination criteria met (e.g., stagnation)} 
            \State \textbf{break}
        \EndIf
    \EndFor
    
    \State \Return The best solution found in $P$
\EndProcedure
\end{algorithmic}
\end{algorithm}

\section{Benchmark Test Instances}
\begin{figure}[h]
    \centering
    \includegraphics[width=1.0\linewidth]{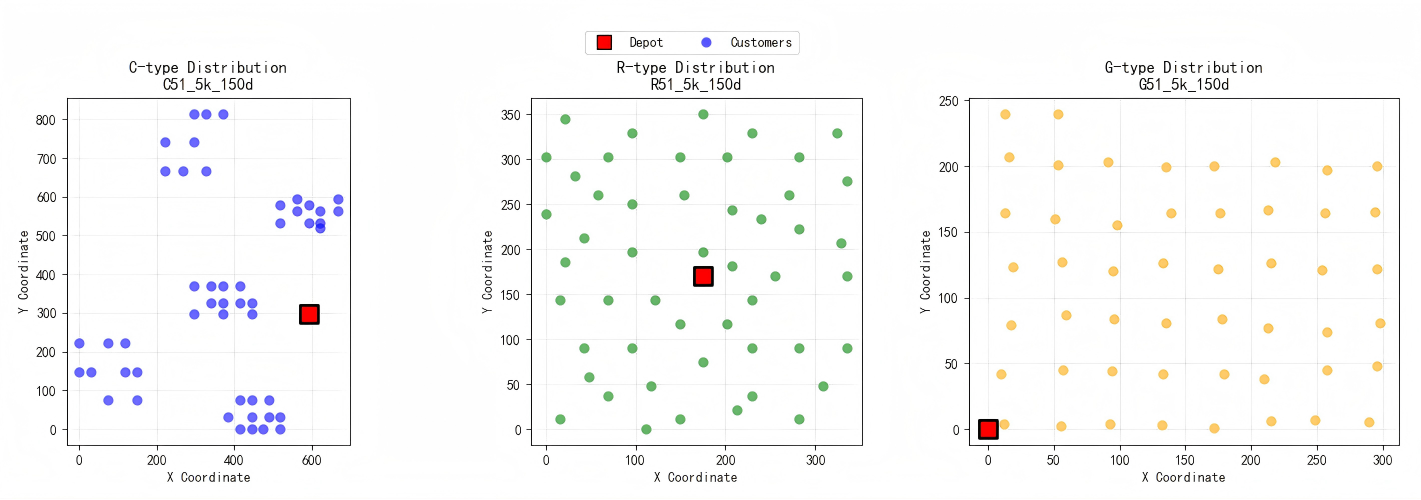}
    \caption{Spatial Distribution Patterns of Different Instance Types: Random (R), Clustered (C), and Mixed (G) distributions showing distinct geographical characteristics of customer locations}
    \label{fig:enter-label}
\end{figure}

The spatial distribution of three types are shown in Figure \ref{fig:enter-label}. The details for the small-scale instances are given in Table \ref{tab:small_instances}, and the medium-scale instances are given in Table \ref{tab:medium_instances}.
\begin{table}[H]
    \centering
    \begin{tabular}{@{}lcccc@{}}
        \toprule
        \textbf{Instance Name} & \textbf{Nodes} & \textbf{K} & \textbf{$\text{D}_\text{max(m)}$} & \textbf{Distribution} \\
        \midrule
        C25\_(2,5)k\_150d & 25 & 2, 5 & 150 & C\\
        R25\_(2,5)k\_150d & 25 & 2, 5 & 150 & R  \\
        G25\_(2,5)k\_150d & 25 & 2, 5 & 150 & G \\
        C/R/G\_50\_5k\_150d & 50 & 5 & 150 & C, R, G \\
        C/R/G\_100\_8k\_150d & 100 & 8 & 150 & C, R, G\\
        G50\_5k\_80d  & 50 & 5 & 80 & G \\
        G50\_5k\_200d & 50 & 5 & 200 & G \\
        G50\_8k\_150d & 50 & 8 & 150 & G \\
        \bottomrule
    \end{tabular}
    \caption{Small Scale STCVRP Instances Dataset (N $\leq$ 100)}
    \label{tab:small_instances}
\end{table}

\begin{table}[H]
    \centering
    \begin{tabular}{@{}lcccc@{}}
        \toprule
        \textbf{Instance Name} & \textbf{Nodes} & \textbf{K} & \textbf{$\text{D}_\text{max(m)}$} & \textbf{Distribution} \\
        \midrule
         G575\_8k\_500d  & 575 & 8  & 500 & G  \\
         G575\_15k\_120d & 575 & 15 & 120 & G \\
         G575\_15k\_500d & 575 & 15 & 500 & G  \\
         G575\_15k\_800d & 575 & 15 & 800 & G \\
         G575\_30k\_500d & 575 & 30 & 500 & G\\
         R575\_15k\_500d & 575 & 15 & 500 & R \\
        \midrule
         R/G\_783\_20k\_500d & 783 & 20 & 500 &  R\\
         G783\_40k\_500d & 783 & 40 & 500 & G \\
         G783\_20k\_(120,800)d & 783 & 20 & 120, 800 & G \\
        \bottomrule
    \end{tabular}
    \caption{Medium Scale STCVRP Instances Dataset (N $>$ 500)}
    \label{tab:medium_instances}
\end{table}

\section{Detailed Exact Algorithm Results}
The detailed computational results of the exact algorithm (Gurobi 12.0.1) on extremely small-scale instances are presented in Table \ref{tab:exact_detailed_results}. These results demonstrate the exponential growth in computational complexity as problem size increases and the extreme sensitivity to constraint parameters.

To simplify the problem for exact solution analysis, all parameters for these extremely small-scale instances were configured with reduced complexity: vehicle speed is set at $1.0$ m/s, service time is set at $8$ seconds, maximum waiting time is set at $6$ seconds, and the inter-task distance is uniformly set at $1$ unit. The only exception is the constraint sensitivity test instance (10-point with doubled distance), where the inter-task distance is set at $2$ units to demonstrate the impact of constraint tightness on computational complexity.
\begin{table}[H]
    \centering
    \begin{tabular}{@{}lccccccc@{}}
        \toprule
        \textbf{Instance} & \textbf{Tasks}  &  \textbf{Time (s)} & \textbf{Nodes Explored} & \textbf{Makespan}\\
        \midrule
        G6\_2k\_30d & 6 &  2.62 & 9,791 & 38.64  \\
        G7\_2k\_30d & 7 & 5.12 & 23528 & 44.07\\
        G8\_2k\_30d & 8 &240.69 & 89,816 & 49.67 \\
        G9\_2k\_30d &9 &1501.43 &995993 & 55.11\\
        G10\_2k\_30d & 10 & 10,912.25 & 9,674,552 & 60.22 \\
        G10*\_2k\_30d & 10* & 32,497.17 & 18,978,649 & 57.77  \\
        \bottomrule
    \end{tabular}
    \caption{Detailed computational results of exact algorithm on extremely small-scale instances}
    \label{tab:exact_detailed_results}
    \vspace{0.2cm}
    \small
    Note: Distance refers to the inter-task distance parameter. The  G10*\_2k\_30d instance has doubled inter-task distances compared to the standard 10-point instance, demonstrating constraint sensitivity.
\end{table}


\section{Performance Analysis of the Genetic Algorithm}
This section will contain convergence curves and performance statistics for the genetic algorithm implementation, demonstrating the algorithm's efficiency and reliability across different instance types and scales.
\begin{figure}[t]
\centering
\begin{subfigure}[b]{0.48\textwidth}
    \centering
    \includegraphics[width=\textwidth]{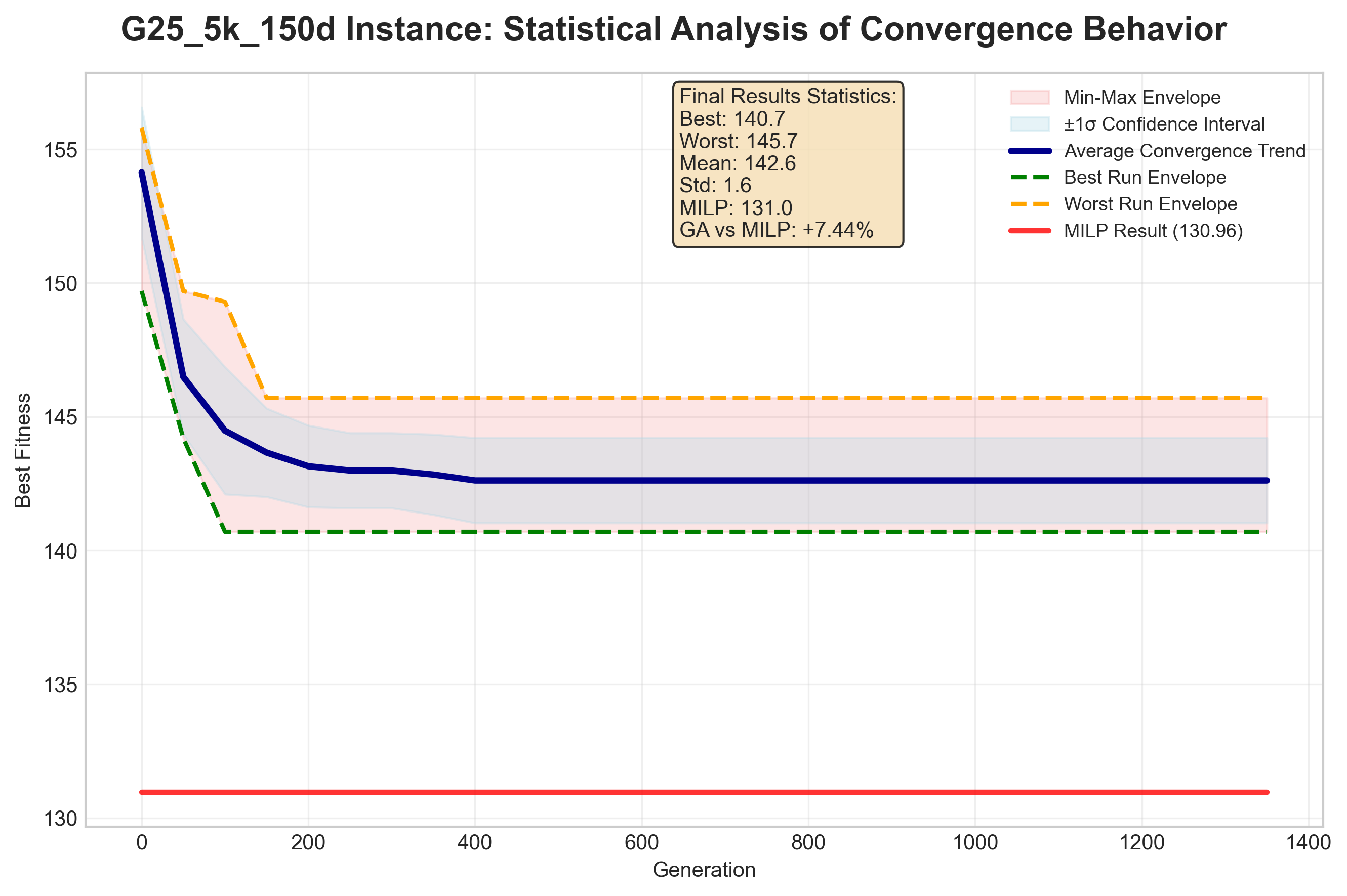}
    \caption{G25\_5k\_150d Instance}
    \label{fig:G25_5k_150d}
\end{subfigure}
\hfill
\begin{subfigure}[b]{0.48\textwidth}
    \centering
    \includegraphics[width=\textwidth]{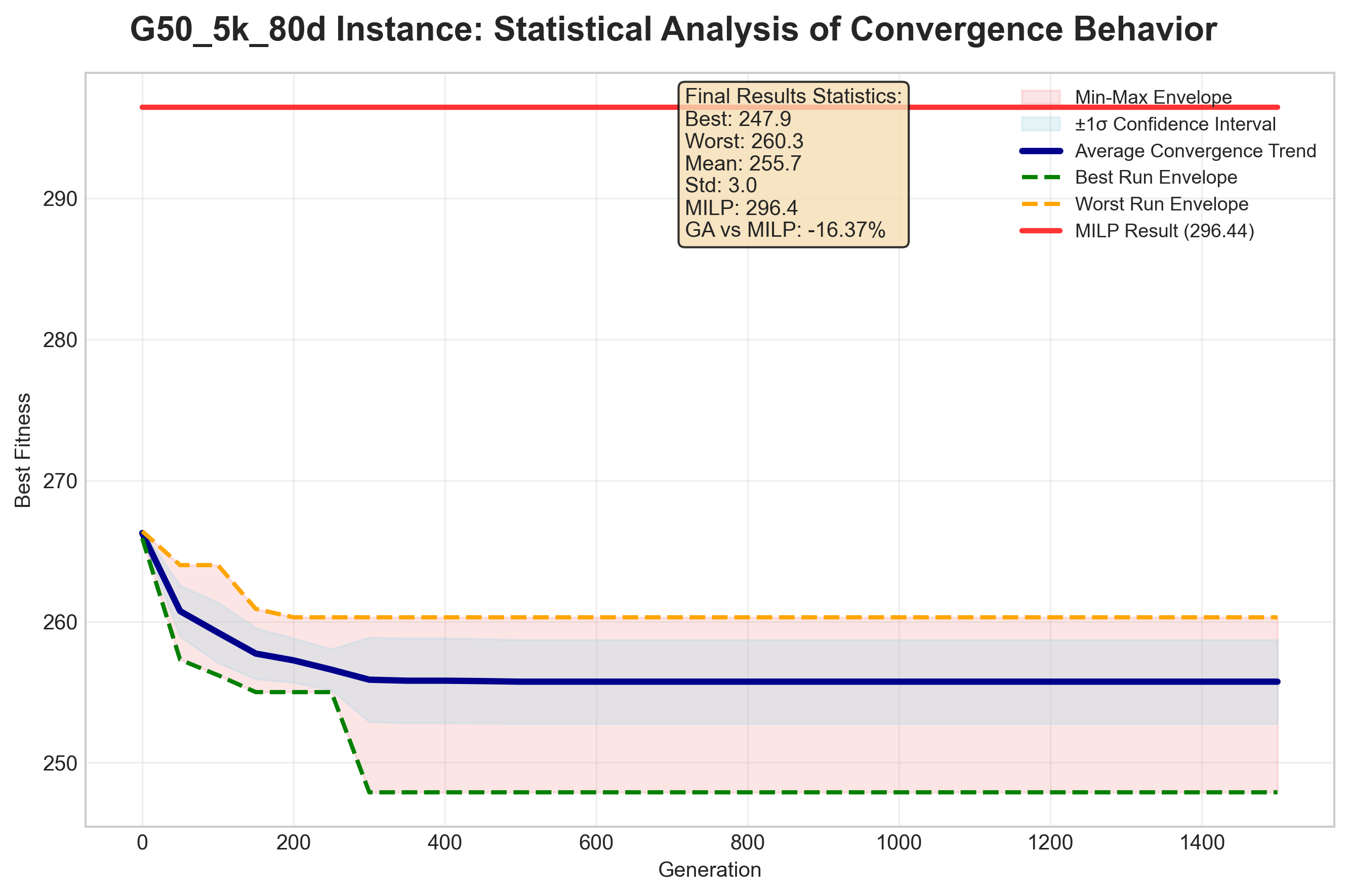}
    \caption{G50\_5k\_80d Instance}
    \label{fig:G50_5k_80d}
\end{subfigure}

\vspace{0.5cm} 

\begin{subfigure}[b]{0.48\textwidth}
    \centering
    \includegraphics[width=\textwidth]{G50_5k_200d.png}
    \caption{G50\_5k\_200d Instance}
    \label{fig:G50_5k_200d}
\end{subfigure}
\hfill
\begin{subfigure}[b]{0.48\textwidth}
    \centering
    \includegraphics[width=\textwidth]{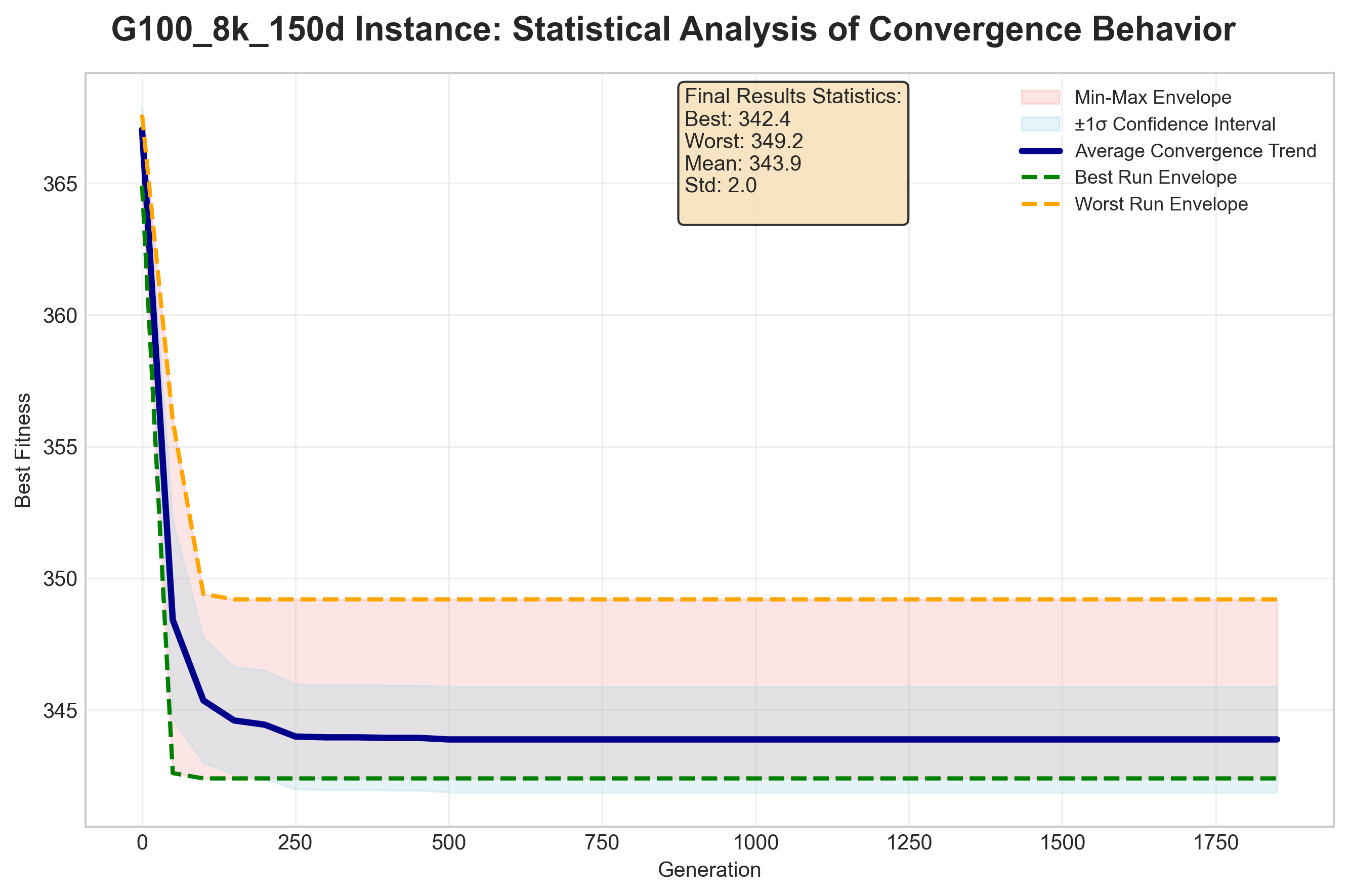}
    \caption{G100\_8k\_150d Instance}
    \label{fig:G100_8k_150d}
\end{subfigure}

\vspace{0.5cm} 


\begin{subfigure}[b]{0.48\textwidth}
    \centering
    \includegraphics[width=\textwidth]{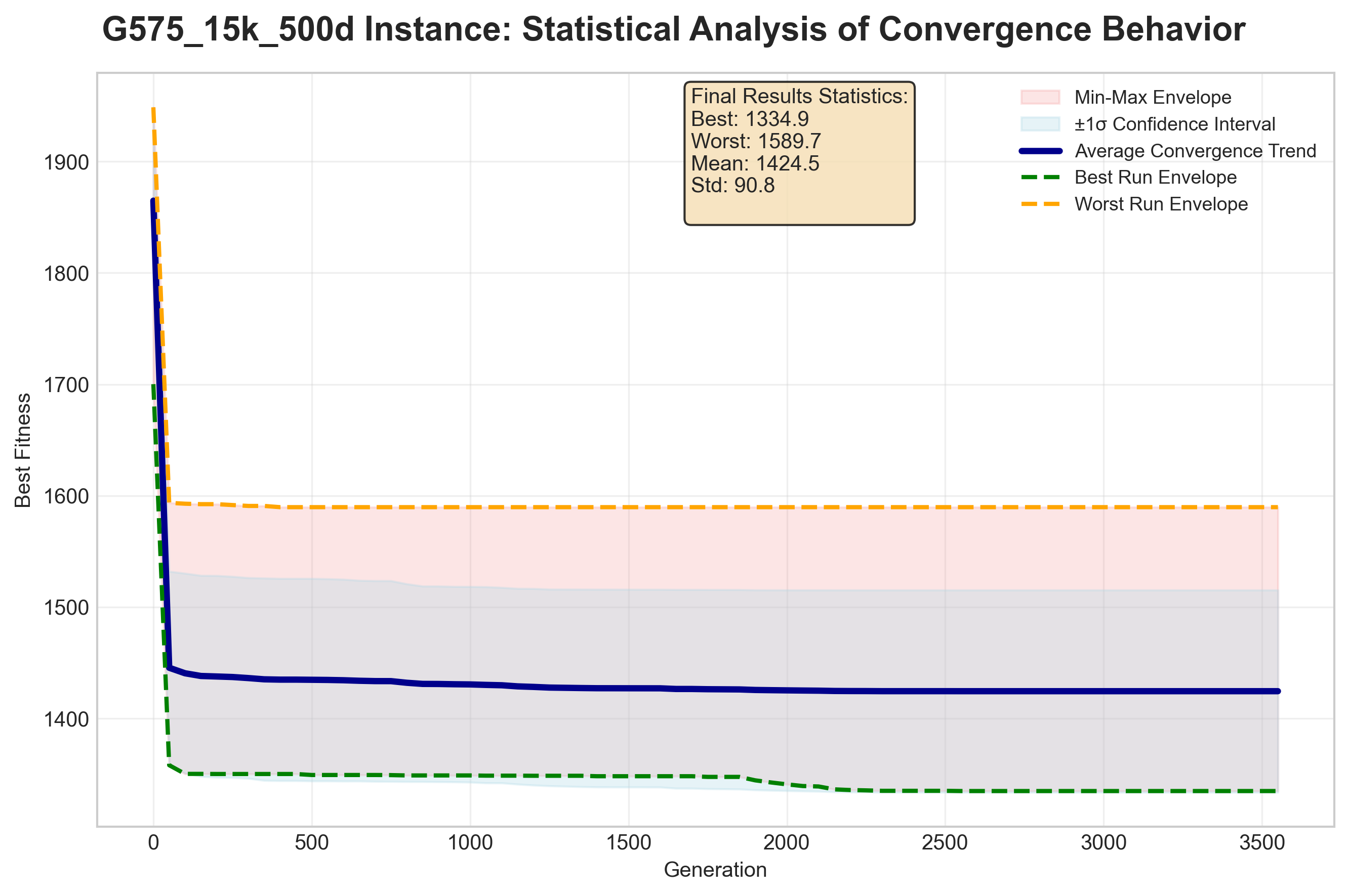}
    \caption{G575\_15k\_500d Instance}
    \label{fig:G575_15k_500d}
\end{subfigure}
\hfill
\begin{subfigure}[b]{0.48\textwidth}
    \centering
    \includegraphics[width=\textwidth]{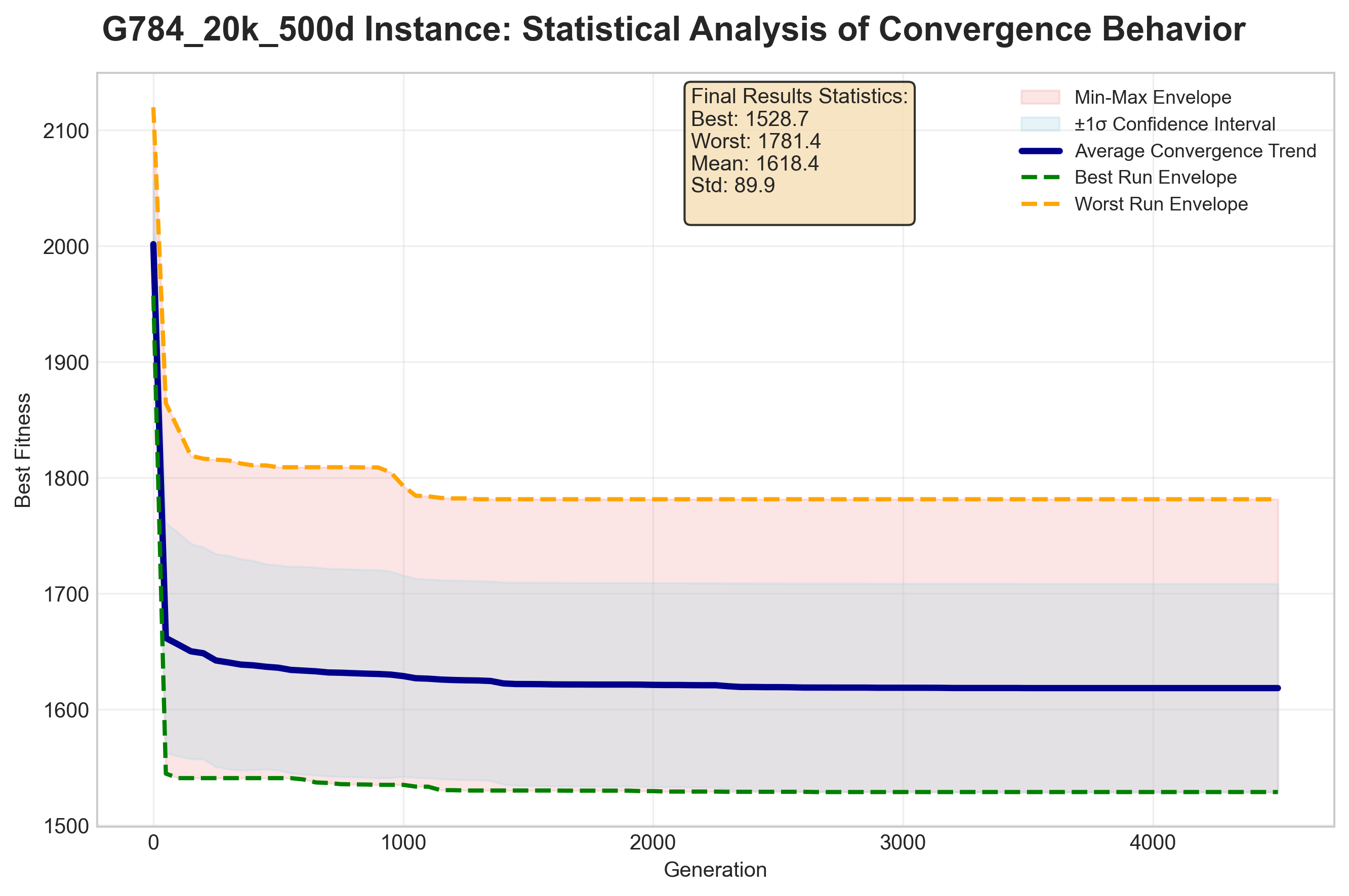}
    \caption{G784\_20k\_500d Instance}
    \label{fig:G784_20k_500d}
\end{subfigure}

\caption{Individual Run Convergence Curves for Different Problem Instances}
\label{fig:all_convergence_curves}
\end{figure}

The iterative convergence patterns, illustrated by representative instances in Figure \ref{fig:all_convergence_curves}, reveal distinct scale-dependent characteristics that provide insights into the algorithm's optimization dynamics. 

Small-scale instances (N $\le$ 100) consistently demonstrate a characteristic two-phase convergence profile. For simpler problems like G25\_5k\_150d, the algorithm exhibits a rapid initial descent, achieving near-optimal solutions within the first 100-200 generations. This is followed by a prolonged fine-tuning phase with minimal further improvements, and the low standard deviation across runs indicates exceptional solution consistency. As the problem size increases within this category to instances like G50\_5k\_150d, this initial convergence period extends, requiring approximately 1000-1500 generations to stabilize, while still maintaining reasonable solution stability.

Medium-scale instances ($N > 500$), such as G575\_15k\_500d and G783\_20k\_500d, exhibit a fundamentally different behavior. These problems require prolonged and continuous optimization phases, often extending beyond 3000 generations, with substantially increased solution variability (standard deviations can be as high as 90.76). Despite this increased variability, all runs demonstrate consistent downward-trending convergence, revealing the algorithm's persistent search capability and an effective balance between exploitation and exploration. 

These convergence curves collectively illustrate that while smaller instances allow for rapid stabilization, larger and more complex instances necessitate extended search processes where the algorithm continues to discover meaningful improvements. This suggests that the GA framework automatically adapts its search intensity to the problem's complexity, effectively trading convergence speed for deeper solution quality exploration in challenging optimization landscapes.

\end{appendix}
\end{document}